\documentclass[10pt]{article}
\usepackage{amsmath}
\usepackage{amsfonts}
\usepackage{amssymb}
\usepackage[english]{babel}
\usepackage{amsthm}

\def\pf{\par\noindent {\bf Proof}~\par\noindent}
\def\qed{~\hfill{$\square$}\pagebreak[1]\par\medskip\par}

\newcommand{\mR}{\mathbb{R}}
\newcommand{\mC}{\mathbb{C}}

\newcommand{\mS}{\mathbb{S}}

\newcommand{\mcP}{\mathcal{P}}
\newcommand{\mcQ}{\mathcal{Q}}
\newcommand{\mcM}{\mathcal{M}}
\newcommand{\mcH}{\mathcal{H}}

\newcommand{\mcX}{\mathcal{X}}
\newcommand{\mcY}{\mathcal{Y}}

\newcommand{\gf}{\mathfrak{f}}
\newcommand{\gfd}{\mathfrak{f}^{\dagger}}
\newcommand{\gsl}{\mathfrak{sl}}

\newcommand{\uX}{\underline{X}}
\newcommand{\ux}{\underline{x}}

\newcommand{\uz}{\underline{z}}
\newcommand{\uzp}{\underline{z}^{\dagger}}
\newcommand{\upz}{\partial_{\uz}}
\newcommand{\wupz}{\widetilde{\upz}}
\newcommand{\upzp}{\partial_{\uz}^{\dagger}}
\newcommand{\wupzp}{\widetilde{\upzp}}

\newcommand{\upzpt}{\widetilde{\upzp}}

\newcommand{\p}{\partial}
\newcommand{\dirac}{\underline{\p}}

\newtheorem{theorem}{Theorem}

\newtheorem{proposition}{Proposition}
\newtheorem{remark}{Remark}
\newtheorem{corollary}{Corollary}

\begin{document}
\title{The Cauchy--Kovalevskaya Extension Theorem in Hermitean Clifford Analysis}
\author{F.\ Brackx$^\ast$, H.\ De Schepper$^\ast$, R.\ L\'{a}vi\v{c}ka$^\ddagger$ \& V.\ Sou\v{c}ek$^\ddagger$}

\date{\small{$^\ast$ Clifford Research Group, Faculty of Engineering, Ghent University\\
Building S22, Galglaan 2, B-9000 Gent, Belgium\\
$^\ddagger$ Mathematical Institute, Faculty of Mathematics and Physics, Charles University\\
Sokolovsk\'a 83, 186 75 Praha, Czech Republic}}

\maketitle

\begin{abstract}
\noindent Hermitean Clifford analysis is a higher dimensional function theory centered around the simultaneous null solutions, called Hermitean monogenic functions, of two Hermitean conjugate complex Dirac operators. As an essential step towards the construction of an orthogonal basis of Hermitean monogenic polynomials, in this paper a Cauchy--Kovalevskaya extension theorem is established for such polynomials. The minimal number of initial polynomials needed to obtain a unique Hermitean monogenic extension is determined, along with the compatibility conditions they have to satisfy. The Cauchy--Kovalevskaya extension principle then allows for a dimensional analysis of the spaces of spherical Hermitean monogenics, i.e.\ homogeneous Hermitean monogenic polynomials. A version of this extension theorem for specific real-analytic functions is also obtained.
\end{abstract}


%%%%%%%%%%%%%%%%%%%%%%%%%%%%%%%%%%%%%%%%%%%%%%%%%%%%%%%%%%%%%%%%%%%%%%%%%%%%%%%%%%%%%%%%%%%%%%%%%%%%%%%%%%%%%%%%%%%%%%%%%%%%%%%%%%%%%%%%%%%%%%%%%%%
\section{Introduction}
%%%%%%%%%%%%%%%%%%%%%%%%%%%%%%%%%%%%%%%%%%%%%%%%%%%%%%%%%%%%%%%%%%%%%%%%%%%%%%%%%%%%%%%%%%%%%%%%%%%%%%%%%%%%%%%%%%%%%%%%%%%%%%%%%%%%%%%%%%%%%%%%%%%%%

The Cauchy--Kovalevskaya theorem (see e.g.\ \cite{cauchy, kowa}) is very well known; for a nice and well--documented historical account on this result we refer to \cite{cooke}. In its most simple setting the theorem reads as follows.

\begin{theorem}
If the functions $g, f_0, \ldots, f_{k-1}$ are analytic in a neighbourhood of the origin, then the initial value problem
\begin{eqnarray*}
\p_t^k h(\ux,t) &=& g(\ux,t,\p_t^i \p_{\ux}^{\alpha} h) \\
\p_t^j h(\ux,0) &=& f_j(\ux), \quad j=0,\ldots,k-1
\end{eqnarray*}
has a unique solution which is analytic in a neighbourhood of the origin, provided that $|\alpha| + i \leq k$.
\end{theorem}
\noindent If the differential operator is chosen to be the Cauchy--Riemann operator, or more explicitly: $\p_t h = i \p_x h$ (with $k=1, |\alpha|=1, i=0$), it follows from this theorem that a holomorphic function in an appropriate region of the complex plane is completely determined by its restriction to the real axis. For a harmonic function though, or more explicitly when $\p^2_t h = -  \p^2_x h$ (with $k=2, |\alpha|=2, i=0$), also the values of its normal derivative on the real axis should be given. In fact, there is a nice and simple construction formula for the holomorphic and harmonic CK--extensions, illustrating the necessity of these restricted values.

\begin{proposition}
If the function $f_0(x)$ is real-analytic in $|x| < a$, then 
\begin{equation}
\label{CK-holo}
F(z) \equiv F(x+it) = \exp{(i t \frac{d}{dx})} \left[f_0(x)\right] = \sum_{k=0}^{\infty} \ \frac{1}{k!} i^k t^k f_0^{(k)}(x)
\end{equation}
is holomorphic in $|z| < a$ and $F(z)|_{t=0} = f_0(x)$. If moreover $f_1(x)$ is real-analytic in $|x| < a$, then 
$$
G(z) \equiv G(x+it) = \sum_{j=0}^{\infty} \ \frac{(-1)^j}{(2j)!} t^{2j} \left(\frac{d}{dx} \right)^{2j} \left[f_0(x)\right] +
\sum_{j=0}^{\infty} \ \frac{(-1)^j}{(2j+1)!} t^{2j+1} \left(\frac{d}{dx} \right)^{2j} \left[f_1(x)\right]
$$
is harmonic in $|z| < a$ and $G(z)|_{t=0} = f_0(x), \frac{\p}{\p t}G(z)|_{t=0} = f_1(x)$.
\end{proposition}

The holomorphic CK--extension principle has been elegantly generalized to higher dimension in the framework of Clifford analysis, which in its most basic form is a higher dimensional generalization of holomorphic function theory in the complex plane, and a refinement of harmonic analysis, see e.g.\ \cite{bds, dss, gilbert, guerleb}. At the heart of this function theory lies the notion of a monogenic function, i.e.\ a Clifford algebra valued null solution of the Dirac operator $\dirac = \sum_{\alpha=1}^m e_{\alpha} \, \p_{X_{\alpha}}$, where $(e_1, \ldots, e_m)$ is an orthonormal basis of $\mR^m$ underlying the construction of the real Clifford algebra $\mR_{0,m}$. We refer to this setting as the Euclidean case, since the fundamental group leaving the Dirac operator $\dirac$ invariant is the orthogonal group $\mbox{O}(m;\mR)$, which is doubly covered by the Pin($m$) group of the Clifford algebra.\\[-3mm]

In the books \cite{rocha,struppa} and the series of papers \cite{sabadini,toulouse,eel,partI,partII,howe,eel2} so--called Hermitean Clifford analysis recently emerged as a refinement of Euclidean Clifford analysis, where the considered functions now take their values in the complex Clifford algebra $\mC_m$ or in complex spinor space. Hermitean Clifford analysis is based on the introduction of an additional datum, a so--called complex structure $J$, inducing an associated Dirac operator $\dirac_J$. Hermitean Clifford analysis then focusses on the simultaneous null solutions of both operators $\dirac$ and $\dirac_J$, called Hermitean monogenic functions. The corresponding function theory is still in full development, see also \cite{hehe, hermwav, mabo, cama, mahi}.\\[-3mm]

The aim of this paper is to establish a CK--extension theorem in Hermitean Clifford analysis and, in particular, to determine the necessary restricted values needed for a unique Hermitean monogenic extension to exist. We confine ourselves to homogeneous polynomials in view of the application of the obtained CK--extension theorem in the construction of orthonormal bases for spaces of Hermitean monogenic homogeneous polynomials using so--called branching rules in group representation theory (see the forthcoming paper \cite{bases}). To make the paper self--contained an introductory section on Clifford analysis is included.

%%%%%%%%%%%%%%%%%%%%%%%%%%%%%%%%%%%%%%%%%%%%%%%%%%%%%%%%%%%%%%%%%%%%%%%%%%%%%%%%%%%%%%%%%%%%%%%%%%%%%%%%%%%%%%%%%%%%%%%%%%%%%%%%%%%%%%%%%%%%%%
\section{Preliminaries of Clifford analysis}
%%%%%%%%%%%%%%%%%%%%%%%%%%%%%%%%%%%%%%%%%%%%%%%%%%%%%%%%%%%%%%%%%%%%%%%%%%%%%%%%%%%%%%%%%%%%%%%%%%%%%%%%%%%%%%%%%%%%%%%%%%%%%%%%%%%%%%%%

For a detailed description of the structure of Clifford algebras we refer to e.g.\ \cite{port}. Here we only recall the necessary basic notions.\\[-3mm]

The real Clifford algebra $\mathbb{R}_{0,m}$ is constructed over the vector space $\mathbb{R}^{0,m}$ endowed with a non--degenerate quadratic form of signature $(0,m)$ and generated by the orthonormal basis $(e_1,\ldots,e_m)$. The non--commutative Clifford or geometric multiplication in $\mathbb{R}_{0,m}$ is governed by the rules 
\begin{equation}\label{multirules}
e_{\alpha} e_{\beta} + e_{\beta} e_{\alpha} = -2 \delta_{\alpha \beta} \ \ , \ \ \alpha,\beta = 1,\ldots ,m
\end{equation}
As a basis for $\mathbb{R}_{0,m}$ one takes for any set $A=\lbrace j_1,\ldots,j_h \rbrace \subset \lbrace 1,\ldots,m \rbrace$ the element $e_A = e_{j_1} \ldots e_{j_h}$, with $1\leq j_1<j_2<\cdots < j_h \leq m$, together with $e_{\emptyset}=1$, the identity element. Any Clifford number $a$ in $\mathbb{R}_{0,m}$ may thus be written as $a = \sum_{A} e_A a_A$, $a_A \in \mathbb{R}$, or still as $a = \sum_{k=0}^m \lbrack a \rbrack_k$, where $\lbrack a \rbrack_k = \sum_{|A|=k} e_A a_A$ is the so--called $k$--vector part of $a$. Euclidean space $\mathbb{R}^{0,m}$ is embedded in $\mathbb{R}_{0,m}$ by identifying $(X_1,\ldots,X_m)$ with the Clifford vector $\uX = \sum_{\alpha=1}^m e_{\alpha}\, X_{\alpha}$, for which it holds that $\uX^2 = - |\uX|^2$. The vector valued first order differential operator $\dirac = \sum_{\alpha=1}^m e_{\alpha}\, \p_{X_{\alpha}}$, called Dirac operator, is the Fourier or Fischer dual of $\uX$. It is this operator which underlies the notion of monogenicity of a function, a notion which is the higher dimensional counterpart of holomorphy in the complex plane. More explicitly, a function $f$ defined and continuously differentiable in an open region $\Omega$ of $\mathbb{R}^{m}$ and taking values in (a subspace of) the Clifford algebra $\mathbb{R}_{0,m}$ is called (left) monogenic in $\Omega$ if $\dirac\lbrack f \rbrack = 0$ in $\Omega$. As the Dirac operator factorizes the Laplacian: $\Delta_m = - \dirac^2$, monogenicity can be regarded as a refinement of harmonicity. The Dirac operator being rotationally invariant, the above framework is usually referred to as Euclidean Clifford analysis. \\[-3mm]

The CK--extension theorem in Euclidean Clifford analysis is a direct generalization to higher dimension of the complex plane case; the monogenic extension is completely determined by its restriction to a real codimension $1$ subspace. It reads as follows (see \cite[p.110]{bds} for the version related to the Cauchy--Riemann operator instead of the Dirac operator used here).
\begin{theorem}
If $\widetilde{f}(X_1,X_2,\ldots,X_{m-1})$ is real-analytic in an open set $\widetilde{\Omega}$ of $\mR^{m-1}$ identified with $\{\uX \in \mR^m: X_m=0\}$, then  there exists an open neigbourhood $\Omega$ of $\widetilde{\Omega}$ in $\mR^{m}$ and a unique monogenic function $f$ in $\Omega$ such that its restriction to $\widetilde{\Omega}$ is precisely $\widetilde{f}$. If moreover $\widetilde{\Omega}$ contains the origin, then in an open neigbourhood of the origin this CK-extension $f = {\rm CK}[\widetilde{f}]$ is given by
\begin{equation}
\label{CKmonog}
f(X_1,X_2,\ldots,X_m) = \exp{\left( X_m e_m \widetilde{\dirac}  \right)} [\widetilde{f}] = \sum_{k=0}^{\infty} \ \frac{1}{k!} X_m^k (e_m \widetilde{ \dirac})^k [\widetilde{f}]
\end{equation}
where $\widetilde{\dirac}$ stands for the restriction of $\dirac$ to $\mR^{m-1}$.
\end{theorem}

Note that the CK-extension operator $\exp{\left( X_m e_m \widetilde{\dirac}  \right)}$ contains the ''new variable'' and derivatives with respect to the ''old'' variables, a typical form which will also be encountered in the Hermitean Clifford analysis case. This CK-operator is the sum of a scalar and a bivector part, since it can be written as
\begin{equation}
\label{CKmonog2}
\exp{\left( X_m e_m \widetilde{\dirac}  \right)} = \sum_{s=0}^{\infty} \ \frac{(-1)^s}{(2s)!} X_m^{2s} \Delta^s + \frac{(-1)^s}{(2s+1)!} X_m^{2s+1} (e_m \widetilde{ \dirac}) \Delta^s
\end{equation}
Note also that when the function $\widetilde{f}$ would happen to be monogenic, i.e. $\widetilde{ \dirac}[\widetilde{f}] =0$, then CK$[\widetilde{f}] = \widetilde{f}$.\\

When allowing for complex constants and taking the dimension to be even: $m=2n$, the generators $(e_1,\ldots, e_{2n})$, still satisfying (\ref{multirules}), produce the complex Clifford algebra $\mathbb{C}_{2n} = \mathbb{R}_{0,2n} \oplus i\, \mathbb{R}_{0,2n}$. Any complex Clifford number $\lambda \in \mathbb{C}_{2n}$ may thus be written as $\lambda = a + i b$, $a,b \in \mathbb{R}_{0,2n}$, leading to the definition of the Hermitean conjugation $\lambda^{\dagger} = (a +i b)^{\dagger} = \overline{a} - i \overline{b}$, where the bar notation stands for the Clifford conjugation in $\mathbb{R}_{0,2n}$, i.e. the main anti--involution for which $\overline{e}_{\alpha} = -e_{\alpha}$, $\alpha=1, \ldots,2n$. This Hermitean conjugation leads to a Hermitean inner product on $\mathbb{C}_{2n}$ given by $(\lambda,\mu) = \lbrack \lambda^{\dagger} \mu \rbrack_0$ and its associated norm $|\lambda| = \sqrt{ \lbrack \lambda^{\dagger} \lambda \rbrack_0} = ( \sum_A |\lambda_A|^2 )^{1/2}$. This is the framework for so--called Hermitean Clifford analysis, where the considered functions will now be defined in open regions of $\mathbb{C}^n \simeq \mathbb{R}^{2n}$ and take their values in (subspaces of) the complex Clifford algebra $\mathbb{C}_{2n}$. \\[-3mm]

An elegant way of introducing this framework consists in considering a complex structure, i.e.\ an $\mbox{SO}(2n;\mR)$--element $J$ for which $J^2=-\mathbf{1}$ (see \cite{partI,partII}). Here, $J$ is chosen to act upon the generators $e_1,\ldots,e_{2n}$ of $\mathbb{C}_{2n}$ as $J[e_j] = -e_{n+j}$ and $J[e_{n+j}] = e_j$, $j=1,\ldots,n$. The projection operators $\frac{1}{2}(\mathbf{1} \pm iJ)$ associated with $J$ produce the main objects of the Hermitean setting by acting upon the corresponding ones in the Euclidean framework. First the so--called Witt basis elements $(\gf_j,\gf_j^{\dagger})^n_{j=1}$ for $\mC_{2n}$ are obtained: 
\begin{eqnarray*}
\gf_j &=& \phantom{-}\frac{1}{2} (\mathbf{1} + iJ) [e_j] \ = \ \phantom{-}\frac{1}{2} (e_j - i \, e_{n+j}), \quad j=1,\ldots,n \\
\gf_j^\dagger &=& -\frac{1}{2} (\mathbf{1} - iJ)[e_j] \ = \ -\frac{1}{2} (e_j + i\, e_{n+j}), \quad j=1,\ldots,n
\end{eqnarray*}
They satisfy the respective Grassmann and duality identities
$$
\gf_j \gf_k + \gf_k \gf_j = \gf_j^{\dagger} \gf_k^{\dagger} + \gf_k^{\dagger} \gf_j^{\dagger} = 0, \quad
\gf_j \gf_k^{\dagger} + \gf_k^{\dagger} \gf_j = \delta_{jk}, \ \ j,k=1,\ldots, n
$$
whence they are isotropic. Next, a vector in $\mR^{0,2n}$ is now denoted by $(x_1, \ldots , x_n, y_1 , \ldots , y_n)$ and identified with the Clifford vector $\uX  =  \sum_{j=1}^n (e_j\, x_j + e_{n+j}\, y_j)$, producing the Hermitean Clifford variables $\uz$ and $\uzp$:
\begin{eqnarray*}
\uz &=& \phantom{-} \frac{1}{2} (\mathbf{1} + iJ) [\uX] \ = \ \sum_{j=1}^n \gf_j\, z_j \\
\uzp &=& -\frac{1}{2} (\mathbf{1} - iJ)[\uX] \ = \ \sum_{j=1}^n \gf_j^{\dagger}\, z_j^c
\end{eqnarray*}
where complex variables $z_j = x_j + i y_j$ have been introduced, with complex conjugates $z_j^c = x_j - i y_j$, $j=1,\ldots,n$. Finally, the Euclidean Dirac operator $\dirac$ gives rise to the Hermitean Dirac operators $\upz$ and $\upzp$:
\begin{eqnarray*}
\upzp &=& \phantom{-}\frac{1}{4} (\mathbf{1} + iJ) [\dirac] \ = \ \sum_{j=1}^n \gf_j\, \p_{z^c_j} \\
\upz &=& -\frac{1}{4} (\mathbf{1} - iJ)[\dirac] \ = \ \sum_{j=1}^n \gf_j^{\dagger}\, \p_{z_j}  
\end{eqnarray*}
involving the Cauchy--Riemann operators $\p_{z_j^c} = \frac{1}{2} (\p_{x_j} + i \p_{y_j})$  and their complex conjugates $\p_{z_j} = \frac{1}{2} (\p_{x_j} - i \p_{y_j})$  in the respective $z_j$--planes, $j=1,\ldots,n$. \\[-3mm]

A continuously differentiable function $g$ in an open region $\Omega$ of $\mathbb{R}^{2n}$ with values in (a subspace of) the complex Clifford algebra $\mathbb{C}_{2n}$ then is called (left) Hermitean monogenic (or h--monogenic) in $\Omega$ if and only if it satisfies in $\Omega$ the system $\upz\, g = 0 = \upzp\, g$, or, equivalently, the system $\dirac \, g = \dirac_J \, g$, with $\dirac_J = J[\dirac]$.\\[-3mm]

Observe that Hermitean vector variables and Dirac operators are isotropic, i.e. $\uz^2 = (\uzp)^2 = 0$ and $\upz^2 = (\upzp)^2 = 0$, whence the Laplacian allows for the decomposition
$$
\Delta_{2n} = 4(\upz \upzp + \upzp \upz) = 4(\upz + \upzp)^2
$$
while also
\begin{displaymath}
(\uz + \uzp)^2 = \uz\, \uzp + \uzp \uz = |\uz|^2 = |\uzp|^2 = |\uX|^2
\end{displaymath}

In the sequel we will consider functions with values in an irreducible representation $\mS$ of the complex Clifford algebra $\mC_m$, usually called  spinor space. To this end first note that, as a vector space, $\mathbb{C}_{2n}$ is isomorphic with the complex Grassmann algebra $\bigwedge_{2n}^{\ast}(\mC^{2n}) = \bigwedge_{2n}^{\ast} ( \gf_1, \gfd_1, \ldots, \gf_n, \gfd_n )$, containing the subspace 
$$
\mbox{$ {\bigwedge}_{n}^{\ast \dagger}(\mC^{n}) =  {\bigwedge}_{n}^{\ast \dagger} ( \gfd_1, \ldots, \gfd_n ) = \bigoplus\limits_{r=0}^{n} ( {\bigwedge}_n^\dagger )^{(r)}$}
$$
where $( \bigwedge_n^\dagger)^{(r)}$ stands for the space of $r$-blades, i.e.
$$
\mbox{$ ( \bigwedge_n^\dagger )^{(r)} = \mbox{span}_{\mC}\left( \gfd_{k_1} \wedge \gfd_{k_2} \wedge \cdots \wedge \gfd_{k_r} : \{ k_1,\ldots,k_r \} \subset \{1,\ldots,n  \}  \right)$}
$$
Spinor space $\mS$ is then realized within the Clifford algebra using a suitable primitive idempotent $I$, say $I = I_1 \ldots I_n$, with $I_j = \gf_j \gf_j^{\dagger}$, $j=1,\ldots,n$. With that choice it holds that $\mS \equiv \mC_{2n} I \cong \mC_n I \cong {\bigwedge}_n^{\ast \dagger} I$, implying that spinor space may decomposed into so--called homogeneous parts as follows:
$$
\mbox{$ \mS = \bigoplus\limits_{r=0}^n \ ( \bigwedge_n^{\dagger} )^{(r)} I$}
$$

%%%%%%%%%%%%%%%%%%%%%%%%%%%%%%%%%%%%%%%%%%%%%%%%%%%%%%%%%%%%%%%%%%%%%%%%%%%%%%%%%%%%%%%%%%%%%%%%%%%%%%%%%%%%%%%%%%%%%%%%%%%%%%%%%%%%%%%%%%%%%%%%%%%%%%%
%%%%%%%%%%%%%%%%%%%%%%%%%%%%%%%%%%%%%%%%%%%%%%%%%%%%%%%%%%%%%%%%%%%%%%%%%%%%%%%%%%%%%%%%%%%%%%%%%%%%%%%%%%%%%%%%%%%%%%%%%%%%%%%%%%%%%%%%%%%%%%

\section{The Hermitean Cauchy-Kovalevskaya extension}

%%%%%%%%%%%%%%%%%%%%%%%%%%%%%%%%%%%%%%%%%%%%%%%%%%%%%%%%%%%%%%%%%%%%%%%%%%%%%%%%%%%%%%%%%%%%%%%%%%%%%%%%%%%%%%%%%%%%%%%%%%%%%%%%%%%%%%%%%%%%%%%

As mentioned in the introduction we will first establish the CK--extension in Hermitean Clifford analysis for polynomials, to which end it clearly suffices to consider polynomials of fixed bidegree in $(\uz,\uzp)$, and taking values in a fixed homogeneous part of spinor space. However, it is clear from the start that, when choosing such a polynomial arbitrarily, it might not have a Hermitean monogenic extension. So we expect the initial polynomials to be subject to additional conditions.\\[-3mm]

In order to investigate this, let $\mcH \mcM_{a,b}^{(r)}$ denote the space of Hermitean monogenic homogeneous polynomials of fixed bidegree $(a,b)$, taking their values in the homogeneous subspace $({\bigwedge}_n^{\dagger})^{(r)} I$ of spinor space. Then we may write any polynomial $M_{a,b} \in \mcH \mcM_{a,b}^{(r)}$ in the form 
$$
M_{a,b} = \sum_{i=0}^a \ \sum_{j=0}^b \frac{z_n^i}{i!} \, \frac{(z_n^{c})^j}{j!} \; p_{a-i,b-j}
$$
where we have singled out the variables $(z_n, z_n^c)$, in view of taking restrictions to $\{z_n=0=z_n^c\}$, identified with $\mC^{n-1}$, later on. The homogeneous polynomials $p_{a-i,b-j}$ thus only contain the variables $(z_1, z_1^c, \ldots, z_{n-1}, z_{n-1}^c)$. Moreover, we split the value space $({\bigwedge}_n^{\dagger})^{(r)} I$ as
$$
\mbox{$({\bigwedge}_n^{\dagger})^{(r)}$} I = \mbox{$({\bigwedge}_{n-1}^{\dagger})^{(r)}$} (\gfd_1, \ldots, \gfd_{n-1}) \, I \ \bigoplus \ \gfd_n \, \mbox{$({\bigwedge}_{n-1}^{\dagger})^{(r-1)}$}  (\gfd_1, \ldots, \gfd_{n-1}) \,  I 
$$
and accordingly decompose the considered functions as $F = F^0 + \gfd_n \, F^1$, where $F^0$ takes values in $({\bigwedge}_{n-1}^{\dagger})^{(r)} (\gfd_1, \ldots, \gfd_{n-1}) \, I$ while $F^1$ takes values in $({\bigwedge}_{n-1}^{\dagger})^{(r-1)} (\gfd_1, \ldots, \gfd_{n-1}) \, I$. The homogeneous polynomial coefficients $p_{a-i,b-j}$ of $M_{a,b}$ are now decomposed as $p_{a-i,b-j} = p_{a-i,b-j}^0 + \gfd_n p_{a-i,b-j}^1$. The resulting components can be organized in a scheme as shown in Figure 1 below.\\[-3mm]

This splitting of values not being possible for $r=0$ nor for $r=n$, we will treat these two exceptional cases separately. So from now on we assume for the general case that $0 < r < n$. \\[-3mm]

\unitlength 0.85cm
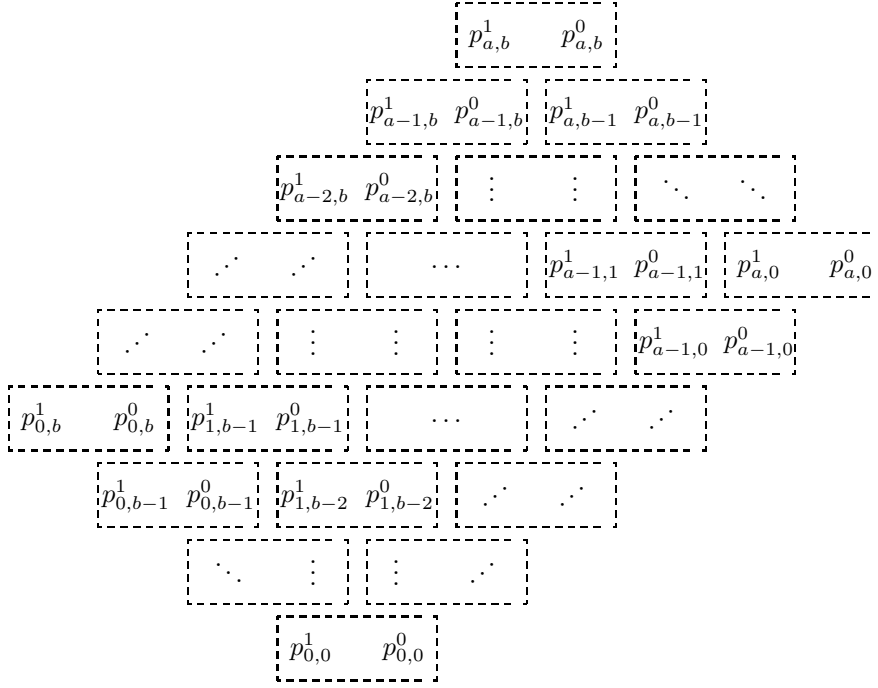
\begin{figure}[h]
\begin{picture}(15,11)(1,0)
\put(5.6,0){\dashbox{0.1}(2.5,1){}}
\put(4.2,1.2){\dashbox{0.1}(2.5,1){}}
\put(7,1.2){\dashbox{0.1}(2.5,1){}}
\put(2.8,2.4){\dashbox{0.1}(2.5,1){}}
\put(5.6,2.4){\dashbox{0.1}(2.5,1){}}
\put(8.4,2.4){\dashbox{0.1}(2.5,1){}}
\put(1.4,3.6){\dashbox{0.1}(2.5,1){}}
\put(4.2,3.6){\dashbox{0.1}(2.5,1){}}
\put(7,3.6){\dashbox{0.1}(2.5,1){}}
\put(9.8,3.6){\dashbox{0.1}(2.5,1){}}
\put(2.8,4.8){\dashbox{0.1}(2.5,1){}}
\put(5.6,4.8){\dashbox{0.1}(2.5,1){}}
\put(8.4,4.8){\dashbox{0.1}(2.5,1){}}
\put(11.2,4.8){\dashbox{0.1}(2.5,1){}}
\put(4.2,6){\dashbox{0.1}(2.5,1){}}
\put(7,6){\dashbox{0.1}(2.5,1){}}
\put(9.8,6){\dashbox{0.1}(2.5,1){}}
\put(12.6,6){\dashbox{0.1}(2.5,1){}}
\put(5.6,7.2){\dashbox{0.1}(2.5,1){}}
\put(8.4,7.2){\dashbox{0.1}(2.5,1){}}
\put(11.2,7.2){\dashbox{0.1}(2.5,1){}}
\put(7,8.4){\dashbox{0.1}(2.5,1){}}
\put(9.8,8.4){\dashbox{0.1}(2.5,1){}}
\put(8.4,9.6){\dashbox{0.1}(2.5,1){}}
\put(8.6,10){$p^1_{a,b}$}
\put(10.05,10){$p^0_{a,b}$}
\put(7.05,8.8){$p^1_{a-1,b}$}
\put(8.38,8.8){$p^0_{a-1,b}$}
\put(9.86,8.8){$p^1_{a,b-1}$}
\put(11.19,8.8){$p^0_{a,b-1}$}
\put(5.65,7.6){$p^1_{a-2,b}$}
\put(6.98,7.6){$p^0_{a-2,b}$}
\put(8.9,7.5){$\vdots$}
\put(10.2,7.5){$\vdots$}
\put(11.6,7.5){$\ddots$}
\put(12.8,7.5){$\ddots$}
\put(8,6.45){$\ldots$}
\put(9.86,6.4){$p^1_{a-1,1}$}
\put(11.19,6.4){$p^0_{a-1,1}$}
\put(12.8,6.4){$p^1_{a,0}$}
\put(14.25,6.4){$p^0_{a,0}$}
\put(6.1,5.1){$\vdots$}
\put(7.4,5.1){$\vdots$}
\put(8.9,5.1){$\vdots$}
\put(10.2,5.1){$\vdots$}
\put(11.26,5.2){$p^1_{a-1,0}$}
\put(12.59,5.2){$p^0_{a-1,0}$}
\put(1.6,4){$p^1_{0,b}$}
\put(3.05,4){$p^0_{0,b}$}
\put(4.26,4){$p^1_{1,b-1}$}
\put(5.59,4){$p^0_{1,b-1}$}
\put(8,4.05){$\ldots$}
\put(2.86,2.8){$p^1_{0,b-1}$}
\put(4.19,2.8){$p^0_{0,b-1}$}
\put(5.66,2.8){$p^1_{1,b-2}$}
\put(6.99,2.8){$p^0_{1,b-2}$}
\put(4.6,1.5){$\ddots$}
%\put(5.8,1.5){$\ddots$}
%\put(7.4,1.45){$\cdot$}
\put(8.6,1.45){$\cdot$}
%\put(7.55,1.58){$\cdot$}
\put(8.75,1.58){$\cdot$}
%\put(7.7,1.71){$\cdot$}
\put(8.9,1.71){$\cdot$}\put(5.8,0.4){$p^1_{0,0}$}
\put(7.25,0.4){$p^0_{0,0}$}
\put(8.8,2.65){$\cdot$}
\put(10,2.65){$\cdot$}
\put(8.95,2.78){$\cdot$}
\put(10.15,2.78){$\cdot$}
\put(9.1,2.91){$\cdot$}
\put(10.3,2.91){$\cdot$}
\put(10.2,3.85){$\cdot$}
\put(11.4,3.85){$\cdot$}
\put(10.35,3.98){$\cdot$}
\put(11.55,3.98){$\cdot$}
\put(10.5,4.11){$\cdot$}
\put(11.7,4.11){$\cdot$}
\put(3.2,5.05){$\cdot$}
\put(4.4,5.05){$\cdot$}
\put(3.35,5.18){$\cdot$}
\put(4.55,5.18){$\cdot$}
\put(3.5,5.31){$\cdot$}
\put(4.7,5.31){$\cdot$}
\put(4.6,6.25){$\cdot$}
\put(5.8,6.25){$\cdot$}
\put(4.75,6.38){$\cdot$}
\put(5.95,6.38){$\cdot$}
\put(4.9,6.51){$\cdot$}
\put(6.1,6.51){$\cdot$}
\put(6.1,1.5){$\vdots$}
\put(7.4,1.5){$\vdots$}
\end{picture}
\caption{Scheme of components of $M_{a,b}$}
\end{figure}

Since it is our intention to investigate the impact of the Hermitean monogenicity of $M_{a,b}$ on its components, we will, in the same order of ideas as above, split the Hermitean variables and Hermitean Dirac operators as
\begin{eqnarray*}
\uz &=& (\gf_1 z_1 + \cdots \gf_{n-1} z_{n-1}) + \gf_n z_n \ = \ \widetilde{\uz} + \gf_n z_n \\
\uzp &=& (\gfd_1 z_1^c + \cdots \gfd_{n-1} z_{n-1}^c) + \gfd_n z_n^c \ = \ \widetilde{\uzp} + \gfd_n z_n^c
\end{eqnarray*}
and
\begin{eqnarray*}
\upz &=& (\gfd_1 \p_{z_1} + \cdots + \gfd_{n-1} \p_{z_{n-1}}) +  \gfd_n \p_{z_n} \ = \ \widetilde{\upz} + \gfd_n \p_{z_n} \\
\upzp &=& (\gf_1 \p_{z_1^c} + \cdots + \gf_{n-1} \p_{z_{n-1}^c}) + \gf_n \p_{z_n^c} \ = \ \widetilde{\upzp} + \gf_n \p_{z_n^c}
\end{eqnarray*}
For $M_{a,b}$ to be Hermitean monogenic, the following conditions should then be satisfied:
\begin{eqnarray*}
\upz M_{a,b} \ = \ \gfd_n \p_{z_n} M_{a,b} + \widetilde{\upz} M_{a,b} & = & 0\\
\upzp M_{a,b} \ = \ \gf_n \p_{z_n^c} M_{a,b} + \widetilde{\upzp} M_{a,b} & = & 0
\end{eqnarray*}
or
\begin{eqnarray*}
\sum_{i=0}^{a-1} \ \sum_{j=0}^{b} \ \frac{z_n^i}{i!} \frac{(z_n^c)^j}{j!} \left( \gfd_n p_{a-i-1,b-j}  \right)  + 
\sum_{i=0}^{a} \ \sum_{j=0}^{b} \ \frac{z_n^i}{i!} \frac{(z_n^c)^j}{j!} \left( \widetilde{\upz} p_{a-i,b-j}\right) & = & 0\\
\sum_{i=0}^{a} \ \sum_{j=0}^{b-1} \ \frac{z_n^i}{i!} \frac{(z_n^c)^j}{j!} \left( \gf_n p_{a-i,b-j-1}  \right)  + 
\sum_{i=0}^{a} \ \sum_{j=0}^{b} \ \frac{z_n^i}{i!} \frac{(z_n^c)^j}{j!} \left( \widetilde{\upzp} p_{a-i,b-j}\right) & = & 0
\end{eqnarray*}
Observe that for $i=a$ and for $j=b$ these conditions, viz $\wupz p_{0,b-j} = 0$, $j=0,\ldots,b$, and $\wupzp p_{a-i,0} = 0$, $i=0,\ldots,a$, are trivially fulfilled. The assumed Hermitean monogenicity of $M_{a,b}$ thus leads to the conditions
\begin{eqnarray}
\label{conditionA}
\gfd_n p_{a-i-1,b-j} + \wupz p_{a-i,b-j} &=& 0 , \quad i=0,\ldots,a-1; \; j=0,\ldots,b\\
\label{conditionB}
\gf_n p_{a-i,b-j-1} + \wupzp p_{a-i,b-j} &=& 0 , \quad i=0,\ldots,a; \; j=0,\ldots,b-1
\end{eqnarray}
Given a minimal number of polynomials on $\mC^{n-1}$, subject to some compatibility conditions, the CK--extension then consists in finding the corresponding unique Hermitean monogenic extension $M_{a,b}$ on $\mC^n$. So we will now determine the necessary and sufficient number of suitable initial polynomials on $\mC^{n-1}$.\\[-3mm]

To this end we make the following observation. From (\ref{conditionA})--(\ref{conditionB}) it follows that 
\begin{eqnarray*}
\wupzp \wupz p_{a-i,b-j} & = & \gfd_n \wupzp p_{a-i-1,b-j} \ = \ - \gfd_n \gf_n p_{a-i-1,b-j-1} \\
\wupz \wupzp p_{a-i,b-j} & = & \gf_n \wupz p_{a-i,b-j-1} \ = \ - \gf_n \gfd_n p_{a-i-1,b-j-1}
\end{eqnarray*}
from which we obtain $(\wupzp \wupz + \wupz \wupzp) p_{a-i,b-j} = - (\gfd_n \gf_n + \gf_n \gfd_n) p_{a-i-1,b-j-1}$, or
\begin{equation}
\label{conditionC}
p_{a-i-1,b-j-1}  =  - \frac{1}{4} \widetilde{\Delta} p_{a-i,b-j}
\end{equation}
whence
$$
p_{a-i-1,b-j-1}^0 = - \frac{1}{4} \widetilde{\Delta} p_{a-i,b-j}^0, \quad p_{a-i-1,b-j-1}^1 = - \frac{1}{4} \widetilde{\Delta} p_{a-i,b-j}^1
$$
since $\widetilde{\Delta}$ is a scalar operator. This means that in the above scheme of Figure 1 each polynomial can be computed directly from the polynomial situated vertically two steps above.\\[-3mm]

Now we analyse conditions (\ref{conditionA}), which we rephrase as
$$
\gfd_n p_{a-i,b-j} + \wupz p_{a-i+1,b-j} = 0 , \quad i=1,\ldots,a; \; j=0,\ldots,b
$$
For $j=0$ and $i=1$ we obtain $\gfd_n p_{a-1,b} + \wupz p_{a,b} = 0$, or
$$
\gfd_n p_{a-1,b}^0 + \gfd_n \gfd_n p_{a-1,b}^1 + \wupz p_{a,b}^0 + \wupz \gfd_n p_{a,b}^1 = 0
$$
The first term and the last term take their values in $\gfd_n \, ({\bigwedge}_{n-1}^{\dagger})^{(r)} (\gfd_1,\ldots,\gfd_{n-1}) \, I$. The second term vanishes due to the isotropy of the Witt basis vectors, while the third term takes its values in $({\bigwedge}_{n-1}^{\dagger})^{(r+1)} (\gfd_1,\ldots,\gfd_{n-1}) \, I$. This leads to a first compatibility condition: $\wupz p_{a,b}^0 = 0$, for $r < n-1$, and a first calculation rule: $p_{a-1,b}^0  = \wupz p_{a,b}^1$. Note that the compatibility condition is trivially satisfied for $r=n-1$ and that the calculation rule implies $\wupz p_{a-1,b}^0 = 0$ due to the isotropy of the Hermitean Dirac operators. For $j=0$ and $i=2$ we obtain $\gfd_n p_{a-2,b} + \wupz p_{a-1,b} = 0$, or
$$
\gfd_n p_{a-2,b}^0 + \gfd_n \gfd_n p_{a-2,b}^1 + \wupz p_{a-1,b}^0 + \wupz \gfd_n p_{a-1,b}^1 = 0
$$
which reduces to $\gfd_n p_{a-2,b}^0 + \wupz \gfd_n p_{a-1,b}^1 = 0$, leading to a second calculation rule: $p_{a-2,b}^0  = \wupz p_{a-1,b}^1$, which implies that also $\wupz p_{a-2,b}^0 = 0$. Proceeding in the same way, still keeping $j=0$ fixed, we find:
$$
p_{a-i,b}^0  = \wupz p_{a-i+1,b}^1 \mbox{\ and \ } \wupz p_{a-i,b}^0 = 0, \quad i=1,\ldots,a
$$
Now taking $j=1$ and $i=1,\ldots,a$ we find the compatibility condition $\wupz p_{a,b-1}^0 = 0$, for $r < n-1$, which again is trivially fulfilled for $r=n-1$, and the calculation rules
$$
p_{a-i,b-1}^0  =  \wupz p_{a-i+1,b-1}^1, \quad i=1,\ldots,a
$$
implying that $\wupz p_{a-i,b-1}^0 = 0$, $i=1,\ldots,a$. \\[-3mm]

Repeating the same reasoning for all values of $j = 0, \ldots, b$ we are lead to
\begin{enumerate}
\item[(i)] the compatibility conditions
$$
\wupz p_{a,b-j}^0 = 0, \quad j=0,\ldots,b, \quad r < n-1
$$
which all are trivially fulfilled for $r=n-1$;
\item[(ii)] the calculation rules
\begin{equation}
\label{conditionE}
p_{a-i,b-j}^0  = \wupz p_{a-i+1,b-j}^1, \quad i=1,\ldots,a; j=0,\ldots,b
\end{equation}
which imply that $\wupz p_{a-i,b-j}^0 = 0$, $i=1,\ldots,a$, $j=0,\ldots,b$.
\end{enumerate}
Note that these compatibility conditions are imposed on the polynomials situated  at the right upper edge of the scheme (see Figure 1). The calculation rules allow for the computation of the $0$-part of each polynomial directly from the $1$-part of the polynomial situated vertically one step up.\\[-3mm]

In a similar way we now analyse the conditions (\ref{conditionB}), which we rephrase as
$$
\gf_n p_{a-i,b-j} + \wupzp p_{a-i,b-j+1} = 0 , \quad i=0,\ldots,a \ ; \ j=1,\ldots,b
$$
For $j=1, i=0$ we obtain $\gf_n p_{a,b-1} + \wupzp p_{a,b} = 0$, or
$$
\gf_n p_{a,b-1}^0 + \gf_n \gfd_n p_{a,b-1}^1 + \wupzp p_{a,b}^0 + \wupzp \gfd_n p_{a,b}^1 = 0
$$
The first term vanishes. The second term can be rewritten as $p_{a,b-1}^1$ which takes its values in $({\bigwedge}_{n-1}^{\dagger})^{(r-1)} (\gfd_1, \ldots, \gfd_{n-1}) \, I$. The third term takes its values in $({\bigwedge}_{n-1}^{\dagger})^{(r-1)} (\gfd_1, \ldots, \gfd_{n-1}) \, I$, and the last term in $\gfd_n \, ({\bigwedge}_{n-1}^{\dagger})^{(r-2)} (\gfd_1, \ldots, \gfd_{n-1}) \, I$. This leads to the compatibility condition $\wupzp p_{a,b}^1  = 0$, for $r > 1$, being trivially satisfied for $r=1$, and the calculation rule
$p_{a,b-1}^1 = - \wupzp p_{a,b}^0$, implying that $\wupzp p_{a,b-1}^1 = 0$. Proceeding in the same way, we finally arrive at
\begin{enumerate}
\item[(iii)] the compatibility conditions 
$$
\wupzp p_{a-i,b}^1 = 0, \quad i=0,\ldots,a, \quad r > 1
$$
which all are trivially fulfilled for $r=1$;
\item[(iv)] the calculation rules
\begin{equation}
\label{conditionG}
p_{a-i,b-j}^1  = - \wupzp p_{a-i,b-j+1}^0 \ , \quad i=0,\ldots,a; \; j=1,\ldots,b
\end{equation}
which imply that $\wupzp p_{a-i,b-j}^1 = 0$, $i=0,\ldots,a$, $j=1,\ldots,b$.
\end{enumerate}
Note that these compatibility conditions are imposed on the polynomials situated  at the left upper edge of the scheme (see Figure 1). The calculation rules allow for the computation of the $1$-part of each polynomial directly from the $0$-part of the polynomial situated vertically one step up.\\[-3mm]

Now (\ref{conditionE}) and (\ref{conditionG}) together yield
\begin{eqnarray*}
p_{a-i,b-j}^1 &=& - \wupzp \wupz p_{a-i+1,b-j+1}^1 \\
p_{a-i,b-j}^0 &=& - \wupz \wupzp p_{a-i+1,b-j+1}^0
\end{eqnarray*}
which is in accordance with (\ref{conditionC}), seen the fact that
$$
\wupz p_{a-i,b-j}^0 = 0 = \wupzp p_{a-i,b-j}^1 I, \quad i=1,\ldots,a; \; j=1,\ldots,b
$$
Summarizing, we have proven the following Hermitean CK--extension theorem.

\begin{theorem}
Given the homogeneous polynomials $p_{a,b-j}^0$, $j=0,\ldots,b$, taking their values in $({\bigwedge}_{n-1}^{\dagger})^{(r)} (\gfd_1, \ldots, \gfd_{n-1}) \, I$ and $p_{a-i,b}^1$, $i=0,\ldots,a$, taking their values in $({\bigwedge}_{n-1}^{\dagger})^{(r-1)} (\gfd_1, \ldots, \gfd_{n-1}) \, I$, satisfying the respective compatibility conditions
\begin{eqnarray*}
&& \wupz p_{a,b}^0 = 0, \ \wupz p_{a,b-1}^0 = 0, \ \ldots, \ \wupz p_{a,0}^0 = 0 \quad (r < n-1) \\
&& \wupzp p_{a,b}^1 = 0, \ \wupzp p_{a-1,b}^1 = 0, \  \ldots, \ \wupzp p_{0,b}^1 = 0 \quad (r > 1)
\end{eqnarray*}
there exists a unique Hermitean monogenic homogeneous polynomial $M_{a,b}$ such that
\begin{itemize}
\item[(i)] $\displaystyle{M_{a,b}|_{\mC^{n-1}} = p_{a,b} = p_{a,b}^0 + \gfd_n p_{a,b}^1}$;
\item[(ii)] $\displaystyle{\frac{\p^j}{\p_{z_n^c}^j} M_{a,b}|_{\mC^{n-1}} = p_{a,b-j} = p_{a,b-j}^0 - \gfd_n \wupzp p_{a,b-j+1}^0}$, $j=1,\ldots,b$;
\item[(iii)] $\displaystyle{\frac{\p^i}{\p_{z_n}^i} M_{a,b}|_{\mC^{n-1}} = p_{a-i,b} =  \wupz p_{a-i+1,b}^1 + \gfd_n p_{a-i,b}^1}$,  $i=1,\ldots,a$.
\end{itemize}
\end{theorem}

\begin{remark}
When $r=0$ all polynomials considered take their values in $({\bigwedge}_{n-1}^{\dagger})^{(0)} \, I$ = $\mbox{span}_{\mC}(1) \, I $, in other words they are scalar polynomials multiplied on the right with the idempotent $I$. It is also known (see \cite{partII}) that in this particular case the notion of Hermitean monogenicity coincides with anti--holomorphy, whence the given polynomials in $\mC^{n-1}$ only depend on the variables $(z_1^c, \ldots, z_n^c)$, so $a$ must be zero. Following the above procedure it turns out that the polynomials $p_{0,b-j}$, $j=0,\ldots,b$, have to be anti--holomorphic, while all other polynomials have to be zero. The corresponding Hermitean monogenic extension is then given by
$$
M_{0,b}^{(0)} = \sum_{j=0}^b \ \frac{{(z_n^c)}^j}{j!} \ p_{0,b-j}
$$
\end{remark}
\begin{remark}
In a similar way, when $r=n$ the polynomials considered take their values in $({\bigwedge}_{n}^{\dagger})^{(n)} \, I$ = $\mbox{span}_{\mC}(\gfd_1\gfd_2\cdots \gfd_n) I$ = $\mbox{span}_{\mC}(\gfd_1\gfd_2\cdots \gfd_n)$, in other words they are scalar polynomials multiplied on the right by $\gfd_1\gfd_2\cdots \gfd_n$. It is also known (see \cite{partII}) that in this particular case the notion of Hermitean monogenicity coincides with holomorphy, whence the given polynomials in $\mC^{n-1}$ only depend on the variables $(z_1, \ldots, z_n)$, so $b$ must be zero. It turns out that the polynomials $p_{a-i,0}$, $i=0,\ldots,a$, have to be holomorphic, while all other polynomials have to be zero. The corresponding Hermitean monogenic extension is then given by
$$
M_{a,0}^{(n)} = \sum_{i=0}^a \ \frac{z_n^i}{i!} \ p_{a-i,0}
$$
\end{remark}

%%%%%%%%%%%%%%%%%%%%%%%%%%%%%%%%%%%%%%%%%%%%%%%%%%%%%%%%%%%%%%%%%%%%%%%%%%%%%%%%%%%%%%%%%%%%%%%%%%%%%%%%%%%%%%%%%%%%%%%%%%%%%%%%%%

\section{Construction formula for the CK--extension}

%%%%%%%%%%%%%%%%%%%%%%%%%%%%%%%%%%%%%%%%%%%%%%%%%%%%%%%%%%%%%%%%%%%%%%%%%%%%%%%%%%%%%%%%%%%%%%%%%%%%%%%%%%%%%%%%%%%%%%%%%%%%%%%%%%%%%

The CK--extension procedure as explained in the preceding section establishes an isomorphism between the space $\mcH \mcM_{a,b}$ of the Hermitean monogenic polynomials of fixed bidegree $(a,b)$, taking their values in a fixed homogeneous subspace of spinor space on the one hand, and the direct sum of the spaces of initially given polynomials subject to compatibility conditions as described in Theorem 3 on the other. Now we want to construct formulae, similar to (\ref{CKmonog}) and (\ref{CKmonog2}), expressing explicitly the Hermitean monogenic extension as the result of the action of an operator on the initially given data. To that end we consider the special case of the foregoing construction where all initially given polynomials but one, are zero. This means that  we execute the CK-extension procedure vertically in one of the columns of the scheme in Figure 1. This leads to the following results.

\begin{corollary}
Given the homogeneous polynomial $p_{a,b}^0$ of fixed bidegree $(a,b)$ with values in $({\bigwedge}_{n-1}^{\dagger})^{(r)} \, I$ and satisfying the compatibility condition $\wupz p_{a,b}^0 = 0$, there exists a unique Hermitean monogenic polynomial $M_{a,b}^0$ such that
\begin{enumerate}
\item[(i)] $\displaystyle{M_{a,b}^0|_{\mC^{n-1}} = p_{a,b}^0}$;
\item[(ii)] $\displaystyle{\p_{z_n^c} M_{a,b}^0|_{\mC^{n-1}} = - \gfd_n \upzpt p_{a,b}^0}$;
\item[(iii)] $\displaystyle{\p_{z_n^c}^j M_{a,b}^0|_{\mC^{n-1}} = 0, \ j = 2, \ldots, b}$;
\item[(iv)] $\displaystyle{\p_{z_n}^i M_{a,b}^0|_{\mC^{n-1}} = 0, \ i = 1, \ldots, a}$.
\end{enumerate}
\label{corol1}
\end{corollary}

\begin{corollary}
Given the homogeneous polynomial $\gfd_n p_{a,b}^1$ of fixed bidegree $(a,b)$ with values in $\gfd_n \, ({\bigwedge}_{n-1}^{\dagger})^{(r-1)} \, I$ and satisfying the compatibility condition $\wupzp p_{a,b}^1 = 0$, there exists a unique Hermitean monogenic polynomial $M_{a,b}^1$ such that:
\begin{enumerate}
\item[(i)] $\displaystyle{M_{a,b}^1|_{\mC^{n-1}} = \gfd_n p_{a,b}^1}$;
\item[(ii)] $\displaystyle{\p_{z_n^c}^j M_{a,b}^1|_{\mC^{n-1}} = 0, \ j = 1, \ldots, b}$;
\item[(iii)] $\displaystyle{\p_{z_n} M_{a,b}^1|_{\mC^{n-1}} = \wupz p_{a,b}^1}$;
\item[(iv)] $\displaystyle{\p_{z_n}^i M_{a,b}^1|_{\mC^{n-1}} = 0, \ i = 2, \ldots, a}$.
\end{enumerate}
\label{corol2}
\end{corollary}

Here, the CK--extension $M_{a,b}^0$ is explicitly given by
$$
M_{a,b}^0 = p_{a,b}^0 + z_n^c \,\gfd_n \, p_{a,b-1}^1 + z_n \, z_n^c \, p_{a-1,b-1}^0 + z_n \, \frac{(z_n^c)^{2}}{2!} \, \gfd_n \, p_{a-1,b-2}^1  + \cdots
$$
or, in view of the calculation rules established in the previous section,
\begin{eqnarray*}
M_{a,b}^0 & = & p_{a,b}^0 + (z_n^c \, \wupzp \, \gfd_n) \, p_{a,b}^0 + (z_n \, \wupz \, \gf_n)(z_n^c \, \wupzp \, \gfd_n) \, p_{a,b}^0 \\
& & + \, \frac{1}{2!}\, (z_n^c \, \wupzp \, \gfd_n)(z_n \, \wupz \, \gf_n)(z_n^c \, \wupzp \, \gfd_n) \, p_{a,b}^0 + \cdots
\end{eqnarray*}
or still, in view of the isotropy of the Hermitean Dirac operators and the compatibility condition for $p_{a,b}^0$,
$$
M_{a,b}^0 = \sum_{k=0}^{\min{(2a+1,2b)}} \frac{1}{\lfloor \frac{k}{2} \rfloor!}\frac{1}{\lfloor \frac{k+1}{2} \rfloor!} \left( z_n \, \wupz \, \gf_n \, +  \, z_n^c \, \wupzp \, \gfd_n \right)^k p_{a,b}^0
$$
which also may be written as
$$
M_{a,b}^0 = \sum_{m=0}^{\min{(a,b)}} \frac{z_n^m}{m!}\frac{(z_n^c)^m}{m!} \left( - \frac{1}{4} \widetilde{\Delta} \right)^m \left[ p_{a,b}^0 \right] \\
+ \sum_{m=0}^{\min{(a,b-1)}} \frac{z_n^m}{m!}\frac{(z_n^c)^{m+1}}{(m+1)!} \left( - \frac{1}{4} \widetilde{\Delta} \right)^m \wupzp \left[ \gfd_n p_{a,b}^0 \right]
$$

Similarly, the CK--extension $M_{a,b}^1$ is explicitly given by
$$
M_{a,b}^1 = \gfd_n \, p_{a,b}^1  + z_n \, p_{a-1,b}^0 + z_n \, z_n^c \, \gfd_n \, p_{a-1,b-1}^1 + \frac{z_n^2}{2!} \, z_n^c\,  p_{a-2,b-1}^0 + \cdots
$$
or, in view of the calculation rules established in the previous section,
\begin{eqnarray*}
M_{a,b}^1 &=& \gfd_n \, p_{a,b}^1 + (z_n \, \wupz \, \gf_n) \, \gfd_n \, p_{a,b}^1 + (z_n^c \, \wupzp \, \gfd_n)(z_n \, \wupz \, \gf_n) \, \gfd_n \, p_{a,b}^1 \\
&& + \, \frac{1}{2!} \, (z_n \, \wupz \, \gf_n)(z_n^c \, \wupzp \, \gfd_n)(z_n \, \wupz \, \gf_n) \gfd_n \, p_{a,b}^1 + \cdots
\end{eqnarray*}
or still, in view of the isotropy of the Hermitean Dirac operators and the compatibility condition for $p_{a,b}^1$,
$$
M_{a,b}^1 = \sum_{k=0}^{\min{(2a,2b+1)}} \frac{1}{\lfloor \frac{k}{2} \rfloor!}\frac{1}{\lfloor \frac{k+1}{2} \rfloor!} \left( z_n \, \wupz \, \gf_n +  z_n^c \, \wupzp \gfd_n \, \right)^k \gfd_n \, p_{a,b}^1
$$
which also may be written as
$$
M_{a,b}^1 = \sum_{m=0}^{\min{(a,b)}} \frac{z_n^m}{m!}\frac{(z_n^c)^m}{m!} \left( - \frac{1}{4} \widetilde{\Delta} \right)^m \left[ \gfd_n p_{a,b}^1 \right] \\
+ \sum_{m=0}^{\min{(a-1,b)}} \frac{z_n^{m+1}}{(m+1)!}\frac{(z_n^c)^{m}}{m!} \left( - \frac{1}{4} \widetilde{\Delta} \right)^m \wupz \left[ p_{a,b}^1 \right]
$$

Finally, the CK--extension $M_{a,b}$ of Theorem 3 is then given by
$$
M_{a,b} = \sum_{j=0}^b \, M_{a,b-j}^0 + \sum_{i=0}^a \, M_{a-i,b}^1
$$
where
$$
M_{a,b-j}^0 = {(z_n^c)}^j \ \sum_{k=0}^{\min{(2a+1,2b-2j)}} \frac{1}{\lfloor \frac{k}{2} \rfloor!}\frac{1}{\lfloor \frac{k+1}{2} + j \rfloor!} \left( z_n \, \wupz \, \gf_n \, +  \, z_n^c \, \wupzp \, \gfd_n \right)^k p_{a,b-j}^0
$$
and
$$
M_{a-i,b}^1 = z_n^i \ \sum_{k=0}^{\min{(2a-2i,2b+1)}} \frac{1}{\lfloor \frac{k}{2} \rfloor!}\frac{1}{\lfloor \frac{k+1}{2} + i \rfloor!} \left( z_n \, \wupz \, \gf_n +  z_n^c \, \wupzp \gfd_n \, \right)^k \gfd_n \, p_{a-i,b}^1
$$
%with $A=2(b-j)$ if $ a \geq b-j$ or $A=2a+1$ if $a<b-j$, and $B=2(a-i)$ if $b \geq a-i$ or $B=2b+1$ if $b<a-i$.

%%%%%%%%%%%%%%%%%%%%%%%%%%%%%%%%%%%%%%%%%%%%%%%%%%%%%%%%%%%%%%%%%%%%%%%%%%%%%%%%%%%%%%%%%%%%%%%%%%%%%%%%%%%%%%%%%%%%%%%%%%%%%%%%%%%%%%

\section{Dimensional analysis}

%%%%%%%%%%%%%%%%%%%%%%%%%%%%%%%%%%%%%%%%%%%%%%%%%%%%%%%%%%%%%%%%%%%%%%%%%%%%%%%%%%%%%%%%%%%%%%%%%%%%%%%%%%%%%%%%%%%%%%%%%%%%%%%%%%%%%%%%

The above CK--extension theorem allows for calculating the dimension of the spaces of Hermitean monogenic polynomials. Indeed, as the CK--extension operator is an isomorphism between the spaces of initially given polynomials subject to  compatibility conditions and the spaces of the corresponding CK--extended Hermitean monogenic polynomials, it will suffice to count the dimension of each of the spaces $\mcX_{a,b,r}$ and $\mcY_{a,b,r}$ consisting of the polynomials $p_{a,b}^0$ and $\gfd_n \, p_{a,b}^1$ respectively. To that end we establish two Fischer decompositions for each of the Hermitean Dirac operators separately.\\[-3mm]

Let $\mcP_{a,b,r}$ be the space of homogeneous polynomials of fixed bidegree $(a,b)$ in the variables $(z_1,z_1^c,\ldots,z_{n-1},z_{n-1}^c)$ with values in $({\bigwedge}_{n-1}^{\dagger})^{(r)} \, I$. The dimension of the value space is $\binom{n-1}{r}$, and the dimension of $\mcP_{a,b,r}$ is given by 
$$
p_{a,b,r} = {\rm dim} (\mcP_{a,b,r}) = {n-1 \choose r} {a+n-2 \choose a} {b+n-2 \choose b} 
$$
Considering the restricted Hermitean Dirac operator $\wupz$ as a linear operator on $\mcP_{a,b,r}$ and denoting $\mcX_{a,b,r} = {\rm Ker} \wupz \cap \mcP_{a,b,r}$, we obtain the following Fischer decomposition for $\wupz$.
\begin{proposition} One has
\begin{equation}
\label{fischerdz}
\mcP_{a,b,r} = \mcX_{a,b,r} \bigoplus \widetilde{\uz} \, \mcX_{a-1,b,r+1}
\end{equation}
\end{proposition}

\pf
Due to the isotropy of the Hermitean Dirac operators, it is clear that Im$\wupz$, being isomorphic with $\mcX_{a,b,r}^{\perp}$, is contained in $\mcX_{a-1,b,r+1}$, whence it suffices to prove that $\mcX_{a-1,b,r+1} \subset {\rm Im}\wupz$. To that end take, for any $\psi \in \mcX_{a-1,b,r+1}$, $\phi = \frac{1}{a+r} \widetilde{\uz} \psi$. Then $\phi \in \mcP_{a,b,r}$ and $\wupz \phi = \psi$. \qed

Now let $\mcQ_{a,b,r}$ be the space of homogeneous polynomials of fixed bidegree $(a,b)$ in the variables $(z_1,z_1^c,\ldots,z_{n-1},z_{n-1}^c)$ with values in $({\bigwedge}_{n-1}^{\dagger})^{(r-1)} \, \gfd_n I$. The dimension of the value space is $\binom{n-1}{r-1}$, and the dimension of $\mcQ_{a,b,r}$ is given by 
$$
q_{a,b,r} = {\rm dim} (\mcQ_{a,b,r}) = {n-1 \choose r-1} {a+n-2 \choose a} {b+n-2 \choose b}
$$
Considering the restricted Hermitean Dirac operator $\wupzp$ as a linear operator on $\mcQ_{a,b,r}$ and denoting $\mcY_{a,b,r} = {\rm Ker} \wupzp \cap \mcQ_{a,b,r}$, we obtain, in a similar way, the following Fischer decomposition for $\wupzp$.
\begin{proposition} One has
\begin{equation}
\label{fischerdzdag}
\mcQ_{a,b,r} = \mcY_{a,b,r} \bigoplus \widetilde{\uzp} \, \mcY_{a,b-1,r-1}
\end{equation}
\end{proposition}

As usual, these Fischer decompositions allow for counting the dimension of the spaces involved. Putting $x_{a,b,r} = {\rm dim} (\mcX_{a,b,r})$ and $y_{a,b,r} = {\rm dim} (\mcY_{a,b,r})$ we deduce from the Fischer decompositions (\ref{fischerdz}) and (\ref{fischerdzdag}) the recurrence relations
\begin{eqnarray*}
p_{a,b,r} &=& x_{a,b,r} + x_{a-1,b,r+1} \\
q_{a,b,r} &=& y_{a,b,r} + y_{a,b-1,r-1}
\end{eqnarray*}
yielding
\begin{eqnarray*}
x_{a,b,r} &=& \frac{r}{a+r} {n-1 \choose r} {a+n-1 \choose a} {b+n-2 \choose b} \\
y_{a,b,r} &=& \frac{r}{b+n-r} {n-1 \choose r} {a+n-2 \choose a} {b+n-1 \choose b}
\end{eqnarray*}
For the dimension of the spaces of Hermitean monogenic polynomials one thus has
$$
m_{a,b}^{(r)} = {\rm dim} (\mcH \mcM_{a,b}^r) = \sum_{j=0}^b \ x_{a,j,r} + \sum_{i=0}^a \ y_{i,b,r}
$$
leading to
\begin{eqnarray*}
m_{a,b}^{(r)} &=& {n-1 \choose r} \left [ \frac{r}{a+r} {a+n-1 \choose a} \sum_{j=0}^b {j+n-2 \choose j} \right . \\
&& \hspace*{3cm} + \left  . \frac{r}{b+n-r} {b+n-1 \choose b} \sum_{i=0}^a {i+n-2 \choose i} \right ]
\end{eqnarray*}
or still
\begin{equation}
\label{dimCK}
m_{a,b}^{(r)} = \frac{(a+b+n) r}{(a+r) (b+n-r)} {n-1 \choose r} {a+n-1 \choose a} {b+n-1 \choose b}
\end{equation}

\unitlength 0.85cm
\begin{figure}[h]
\begin{picture}(13.5,13.5)(-0.2,0)
\put(0,0){\line(1,0){6}}
\put(0,1.5){\line(1,0){6}}
\put(0,3){\line(1,0){6}}
\put(0,4.5){\line(1,0){6}}
\put(0,6){\line(1,0){7.5}}
\put(0,7.5){\line(1,0){7.5}}
\put(0,9){\line(1,0){7.5}}
\put(0,10.5){\line(1,0){7.5}}
\put(0,12){\line(1,0){13.5}}
\put(0,13.5){\line(1,0){13.5}}
\put(0,0){\line(0,1){13.5}}
\put(1.5,0){\line(0,1){13.5}}
\put(3,0){\line(0,1){13.5}}
\put(4.5,0){\line(0,1){13.5}}
\put(6,0){\line(0,1){13.5}}
\put(7.5,6){\line(0,1){7.5}}
\put(9,12){\line(0,1){1.5}}
\put(10.5,12){\line(0,1){1.5}}
\put(12,12){\line(0,1){1.5}}
\put(13.5,12){\line(0,1){1.5}}
\put(0.65,0.65){\footnotesize{$b$}}
\put(2.05,0.7){\footnotesize{$\ldots$}}
\put(3.65,0.65){\footnotesize{$2$}}
\put(5.15,0.65){\footnotesize{$1$}}
\put(0.35,2.1){\footnotesize{$b+1$}}
\put(2.05,2.2){\footnotesize{$\ldots$}}
\put(3.65,2.1){\footnotesize{$3$}}
\put(5.15,2.1){\footnotesize{$2$}}
\put(0.70,3.6){\footnotesize{$\vdots$}}
\put(2.20,3.6){\footnotesize{$\vdots$}}
\put(3.70,3.6){\footnotesize{$\vdots$}}
\put(5.20,3.6){\footnotesize{$\vdots$}}
\put(0.08,5.30){\footnotesize{$b+n$}}
\put(0.33,4.98){\footnotesize{$-r-2$}}
\put(2.05,5.2){\footnotesize{$\ldots$}}
\put(3.35,5.15){\footnotesize{$n-r$}}
\put(4.58,5.15){\footnotesize{$n-r-1$}}
\put(0.08,6.65){\footnotesize{$b+n-r$}}
\put(2.05,6.7){\footnotesize{$\ldots$}}
\put(3.08,6.65){\footnotesize{$n-r+2$}}
\put(4.58,6.65){\footnotesize{$n-r+1$}}
\put(6.65,6.65){\footnotesize{$1$}}
\put(0.1,8.30){\footnotesize{$b+n$}}
\put(0.35,7.96){\footnotesize{$-r+1$}}
\put(2.05,8.2){\footnotesize{$\ldots$}}
\put(3.08,8.15){\footnotesize{$n-r+3$}}
\put(4.58,8.15){\footnotesize{$n-r+2$}}
\put(6.65,8.15){\footnotesize{$2$}}
\put(0.70,9.6){\footnotesize{$\vdots$}}
\put(2.20,9.6){\footnotesize{$\vdots$}}
\put(3.70,9.6){\footnotesize{$\vdots$}}
\put(5.20,9.6){\footnotesize{$\vdots$}}
\put(6.70,9.6){\footnotesize{$\vdots$}}
\put(0.07,11.15){\footnotesize{$b+n-2$}}
\put(2.05,11.2){\footnotesize{$\ldots$}}
\put(3.65,11.15){\footnotesize{$n$}}
\put(4.87,11.15){\footnotesize{$n-1$}}
\put(6.37,11.15){\footnotesize{$r-1$}}
\put(0.1,12.77){\footnotesize{$a+b$}}
\put(0.35,12.43){\footnotesize{$+n-1$}}
\put(3.09,12.65){\footnotesize{$a+n+1$}}
\put(4.89,12.65){\footnotesize{$a+n$}}
\put(6.39,12.65){\footnotesize{$a+r$}}
\put(2.05,12.7){\footnotesize{$\ldots$}}
\put(8.15,12.65){\footnotesize{$a$}}
\put(9.55,12.7){\footnotesize{$\ldots$}}
\put(11.15,12.65){\footnotesize{$2$}}
\put(12.65,12.65){\footnotesize{$1$}}
\end{picture}
\caption{Ferrer diagram with hook numbers for $\mcH\mcM_{a,b}^{(r)}$}
\end{figure}

We will now compare this expression with the dimension formula for spherical Hermitean monogenics obtained in \cite{fischer}.
There the dimension was established following two approaches, one of them involving the Fischer decomposition of harmonic homogeneous polynomials in terms of Hermitean monogenic ones, the other being based on the Weyl dimension formula (see \cite[p.382]{good}). In the second approach the space  $\mcM_{a,b}^{(r)}$ is considered as a $\mathfrak{u}(n)$-module with highest weights 
$$
\lambda = [\lambda_1, \lambda_2, \ldots, \lambda_r, \lambda_{r+1}, \ldots, \lambda_{n-1}] = [a+b+1, b+1, \ldots, b+1, b, \ldots, b]
$$
(see \cite{damdeef}), where the last $b+1$ appears at the $r$-th place. According to the Weyl dimension formula, its dimension is then given by
$$
\frac{(\lambda_1+n-1)!}{(n-1)!} \frac{(\lambda_{2}+n-2)!}{(n-2)!} \cdots \frac{(\lambda_{n-1}+1)!}{1!} \frac{1}{\Pi_{i,j \in \lambda} h_{i,j}}
$$
where the $h_{i,j}$ are the so--called hook numbers shown in the Ferrer diagram above (see Figure 2). \\[-3mm]

In \cite{fischer} this lead for $0 < r < n$ to 
$$
m_{a,b}^{(r)} = \frac{a+b+n}{a+r} \binom{b+n-r-1}{b} \binom{b+n-1}{r-1} \binom{a+n-1}{a}
$$
which indeed exactly coincides with (\ref{dimCK}).

%%%%%%%%%%%%%%%%%%%%%%%%%%%%%%%%%%%%%%%%%%%%%%%%%%%%%%%%%%%%%%%%%%%%%%%%%%%%%%%%%%%%%%%%%%%%%%%%%%%%%%%%%%%%%%%%%%%%%%%%%%%%%%%%%%%%%%%%%%%

\section{CK--extension of a real--analytic function}

%%%%%%%%%%%%%%%%%%%%%%%%%%%%%%%%%%%%%%%%%%%%%%%%%%%%%%%%%%%%%%%%%%%%%%%%%%%%%%%%%%%%%%%%%%%%%%%%%%%%%%%%%%%%%%%%%%%%%%%%%%%%%%%%%%%%%%%%%%

Although the primary aim of this paper was to study the Cauchy-Kovalevskaya extension of homogeneous polynomials, in view of the construction of an orthogonal basis of spaces of Hermitean monogenic homogeneous polynomials, as explained in the introduction, we devote a small section to the CK--extension of a real--analytic function, which is now readily obtained from the results in Section 4.\\

So consider a real--analytic function $\widetilde{F}(x_1, \ldots, x_{n-1}, y_1, \ldots, y_{n-1})$ in the neighbourhood of the origin in $\mR^{2n-2}$ with values in the homogeneous subspace $({\bigwedge}_n^{\dagger})^{(r)} I$ of spinor space $(0 < r < n)$. This function $\widetilde{F}$ may be rewritten as a function of the variables $(z_1, \ldots, z_{n-1}, z_1^c, \ldots, z_{n-1}^c)$, and, as before, it may be decomposed as
$$
\widetilde{F} = \widetilde{F}^0 + \gfd_n \, \widetilde{F}^1
$$
where $\widetilde{F}^0$ takes values in $({\bigwedge}_{n-1}^{\dagger})^{(r)} (\gfd_1, \ldots, \gfd_{n-1}) \, I$ while $\widetilde{F}^1$ takes values in $({\bigwedge}_{n-1}^{\dagger})^{(r-1)} (\gfd_1, \ldots, \gfd_{n-1}) \, I$.\\
As the function $\widetilde{F}$ is assumed to be real--analytic in the neighbourhood of the origin, it may be developed into a convergent multiple power series, first in the variables $(x_1, \ldots, x_{n-1}, y_1, \ldots, y_{n-1})$, then in the variables $(z_1, \ldots, z_{n-1}, z_1^c, \ldots, z_{n-1}^c)$, ending up with
$$
\widetilde{F}(\widetilde{\uz},\widetilde{\uzp}) = \sum_{a=0}^{\infty} \ \sum_{b=0}^{\infty} \ p^0_{a,b} + \gfd_n p^1_{a,b}
$$
with $(\alpha=0,1)$
$$
p^{\alpha}_{a,b} = \sum_{|(a_1, \ldots, a_{n-1})| = a} \ \sum_{|(b_1, \ldots, b_{n-1})| = b} \ z_1^{a_1} \cdots z_{n-1}^{a_{n-1}}(z_1^c)^{b_1} \cdots (z_{n-1}^c)^{b_{n-1}} C^{\alpha}_{(a_1, \ldots, a_{n-1}),(b_1, \ldots, b_{n-1})}
$$
the $C^{\alpha}_{(a_1, \ldots, a_{n-1}),(b_1, \ldots, b_{n-1})}$ being appropriate spinor valued constants.\\

In order for the function $\widetilde{F}$ to have a Hermitean monogenic extension to a neighbourhood of the origin in $\mR^{2n}$, it is clear that its components $\widetilde{F}^0$ and $\widetilde{F}^1$ should simultaneously satisfy the conditions
$$
\wupz \widetilde{F}^0 = 0 \quad , \quad \wupzp \widetilde{F}^1 = 0
$$
The homogeneous polynomials $p^0_{a,b}$ and $p^1_{a,b}$ appearing in the above multiple power series expansion, then precisely satisfy the compatibility conditions of Corollaries 1 and 2, in order to possess a Hermitean monogenic polynomial extension. This leads to the multiple power series 
$$
\sum_{a=0}^{\infty} \ \sum_{b=0}^{\infty} \ CK\left[ p^0_{a,b} \right] + CK\left[ \gfd_n \, p^1_{a,b} \right] 
$$
$$
= \sum_{a=0}^{\infty} \ \sum_{b=0}^{\infty} \ \sum_k \ \frac{1}{\lfloor \frac{k}{2} \rfloor!}\frac{1}{\lfloor \frac{k+1}{2} \rfloor!} \left( z_n \, \wupz \, \gf_n +  z_n^c \, \wupzp \gfd_n \, \right)^k \left[ p^0_{a,b} + \gfd_n \, p_{a,b}^1 \right]
$$
which clearly converges to a Hermitean monogenic function $F(\uz,\uzp)$ in the neighbourhood of the origin in $\mR^{2n}$. Obviously $F(\uz,\uzp)$ is an extension of the given real--analytic function $\widetilde{F}(\widetilde{\uz},\widetilde{\uzp})$.

%%%%%%%%%%%%%%%%%%%%%%%%%%%%%%%%%%%%%%%%%%%%%%%%%%%%%%%%%%%%%%%%%%%%%%%%%%%%%%%%%%%%%%%%%%%%%%%%%%%%%%%%%%%%%%%%%%%%%%%%%%%%%%%%%%%%%

\section*{Acknowledgements}

%%%%%%%%%%%%%%%%%%%%%%%%%%%%%%%%%%%%%%%%%%%%%%%%%%%%%%%%%%%%%%%%%%%%%%%%%%%%%%%%%%%%%%%%%%%%%%%%%%%%%%%%%%%%%%%%%%%%%%%%%%%%%%%%%%%%

R. L\'{a}vi\v{c}ka and V. Sou\v{c}ek acknowledge support by the institutional grant MSM 0021620839 and by grant GA CR 201/08/0397.

%%%%%%%%%%%%%%%%%%%%%%%%%%%%%%%%%%%%%%%%%%%%%%%%%%%%%%%%%%%%%%%%%%%%%%%%%%%%%%%%%%%%%%%%%%%%%%%%%%%%%%%%%%%%%%%%%%%%%%%%%%%%%%%%%%

%%%%%%%%%%%%%%%%%%%%%%%%%%%%%%%%%%%%%%%%%%%%%%%%%%%%%%%%%%%%%%%%%%%%%%%%%%%%%%%%%%%%%%%%%%%%%%%%%%%%%%%%%%%%%%%%%%%%%%

\end{document}